\documentclass[12pt]{article}
\usepackage{amssymb,amsmath}
\usepackage{graphicx}
\usepackage{amsfonts}
\usepackage{float}
\usepackage{flafter}
\usepackage{color}

\labelwidth\leftmargini\advance\labelwidth-\labelsep

\def\R{\relax\ifmmode I\!\!R\else$I\!\!R$\fi}

\def\Z{\relax\ifmmode Z\!\!\!Z\else$Z\!\!\!Z$\fi}

\def\C{\relax\ifmmode C\!\!\!\!I\else$C\!\!\!\!I$\fi}

\def\K{\relax\ifmmode I\!\!K\else$I\!\!K$\fi}

\def\N{\relax\ifmmode I\!\!N\else$I\!\!N$\fi}

\newcounter{defcounter}[section]
{\vspace{0.1cm}\begin{sloppypar}\noindent\stepcounter{defcounter}{\bfseries
Definition
      \thesection.\thedefcounter}}%
{\end{sloppypar}\vspace{0.1cm}}

\newtheorem{lemma}{Lemma}[section]
\newtheorem{theorem}{Theorem}[section]

\newtheorem{proposition}{Proposition}[section]

\newcommand{\proof}{{\bf Proof.} }

\newcommand{\qed}{\hfill $\square$}

\begin{document}
\thispagestyle{empty}
\begin{center}
{\Large {\bf On the universality of Somos' constant}}
\end{center}
\begin{center}J\"org Neunh\"auserer\\
Technical University of Braunschweig \\
joerg.neunhaeuserer@web.de
\end{center}
\begin{center}
\begin{abstract}
We show that Somos' constant is universal in sense that is similar to the universality of the Khinchin constant. In addition we introduce generalized Somos' constants, which are universal in a similar sense.
\end{abstract}
\end{center}
\section{Introduction and main result}
Let us first recall the Khinchin constant
\[ K=\prod_{i=1}^\infty \left(1+\frac{1}{i(i+2)}\right)^{\log_2 i} = 2.6854520010\dots~~.\]
By the famous theorem of Khinchin \cite{[KI]} this constant is universal in the following sense: For almost all real numbers $x$ the geometric mean of the entries of the continued fractions of $x$ converges to $K$. We consider here Somos' constant
\[ \sigma=\prod_{i=1}^\infty \sqrt[2^i]{i}= 1.6616879496\dots~~,\]
which first appeared in \cite{[SO]} in the context of the quadratic recurrence $g_{n}=ng_{n-1}^{2}$, see also page 446 of \cite{[FI]}.
In the recent past this constant raised some attention, see for instance \cite{[GS],[M],[SL]}. We will show that the Somos' constant is universal in a sense that is similar to the universality of the Khinchin constant. In \cite{[NE]} we represent real numbers $x\in(0,1]$ in the form
\[ x=\langle n_{1},n_{2},n_{3},\dots\rangle:=\sum_{k=1}^{\infty}2^{-(n_{1}+n_{2}+\dots+n_{k})}\]
with $n_{k}\in\mathbb{N}$ and show that the representation is unique. Replacing the continued fraction representation by this representation, we obtain the universality of Somos' constant.
\begin{theorem}
For almost all $x=\langle n_{1},n_{2},n_{3},\dots\rangle\in(0,1]$ we have
\[ \lim_{i\to\infty} \sqrt[i]{n_{1}n_{2}\dots n_{i}}=\sigma.\]
\end{theorem}
We will prove this theorem in the next section. In the last section we will introduce generalized Somos' constants, which are universal with respect to a modification of the representation used here.
\section{Proof of the main result}
Consider the map $T:(0,1]\to(0,1]$, given by $T(x)=2^{i}x-1$ for $x\in (1/2^{i},1/2^{i-1}]$ and $i\in\mathbb{N}$. The relation of this transformation to the expansion of real numbers, defined in the last section is given by
\begin{lemma} Let $x=\langle n_{1},n_{2},n_{3},\dots\rangle\in(0,1]$.
For all $k\in\mathbb{N}$ we have $T^{k-1}(x)\in (1/2^{i},1/2^{i-1}]$ if and only if $n_{k}=i$.
\end{lemma}
\proof Obviously $T(\langle n_{1},n_{2},n_{3},\dots\rangle)=\langle n_{2},n_{3},n_{4},\dots\rangle$. Since $x\in (1/2^{i},1/2^{i-1}]$ if and only if $n_{1}=i$ the result follows immediately. \qed\\~\\
To apply Birkhoff's ergodic theorem we prove:
\begin{proposition}
The Lebesgue measure $\mathfrak{L}$ is ergodic with respect to $T$.
\end{proposition}
\proof
For an open interval $(a,b)\subseteq [0,1]$ we have
\[ \mathfrak{L}(T^{-1}((a,b)))=\mathfrak{L}\left(\bigcup_{i=1}^{\infty}(a/2^{i}+1/2^{i},b/2^{i}+1/2^{i})\right)\]
\[ = \sum_{i=1}^{\infty}2^{-k}\mathfrak{L}\left((a/2^{i}+1/2^{i},b/2^{i}+1/2^{i})\right)=\sum_{i=1}^{\infty}2^{-i}(b-a)=b-a=\mathfrak{L}((a,b)).\]
Hence $\mathfrak{L}(T^{-1}(B))=\mathfrak{L}(B)$ for all Borel sets $B\subseteq (0,1]$, which means that $\mathfrak{L}$ is invariant under $T$.
Let $B$ be a Borel set with $\mathfrak{L}(B)<1$, which is invariant under $T$; that is $T(B)=B$.
Note that for all $k\in\mathbb{N}$ the intervals of the form
\[ I_{m_{1},\dots,m_{k}}=\{\langle n_{1},n_{2},n_{3},\dots\rangle|n_{i}=m_{i}\mbox{ for }i=1,\dots,k\}\]
build a partition of $(0,1]$, where the length of the partition elements is bounded by $1/2^{k}$. By Lebesgue's density theorem for every $\epsilon>0$ there is an interval $I= I_{m_{1},\dots,m_{k}}$ such that $\mathfrak{L}(I\backslash B)\ge (1-\epsilon)\mathfrak{L}(I)$. Since $T^{k}(I)=(0,1]$ we have
\[ \mathfrak{L}((0,1]\backslash B)\ge \mathfrak{L}(T^{k}(I\backslash B))\ge (1-\epsilon)\mathfrak{L}(T^{k}(I))= 1-\epsilon.\]
Hence $\mathfrak{L}(B)=0$. This proves that $\mu$ is ergodic.
\qed.\\~\\
Now we are prepared to prove Theorem 1.1. Let $f(x)=\sum_{i=1}^{\infty}\log(i)\chi_{(1/2^{i},1/2^{i-1}]}(x)$, where $\chi$ is the characteristic function.
By lemma 2.1  we have $f(T^{k-1}(x))=\log(n_{k})$ for $x=\langle n_{1},n_{2},n_{3},\dots\rangle$.
Applying Birkhoff's ergodic theorem to $T$ with the $L^1$-function $f$, we obtain
\[ \lim_{i\to\infty}  \frac{1}{i}\sum_{k=1}^{i}\log(n_{k})=\lim_{i\to\infty}\frac{1}{i}\sum_{k=1}^{i}f(T^{k-1}(x))=\int_{0}^{1} f(x)dx \]
\[ = \sum_{i=1}^{\infty}\log(i)2^{-i}\]
for almost all $x=\langle n_{1},n_{2},n_{3},\dots\rangle\in(0,1]$. Taking the exponential gives the result.
\section{A generalisation} Let $b\ge 2$ be an integer.
It is easy to show that a real numbers $x\in(0,1]$ has a unique representation in the form
\[ x=\langle n_{1},n_{2},n_{3},\dots\rangle_{b}:=(b-1)\sum_{k=1}^{\infty}b^{-(n_{1}+n_{2}+\dots+n_{k})}\]
with $n_{k}\in\mathbb{N}$, the argument can be found in $\cite{[NE]}$. Now consider the map $T_{b}:(0,1]\to(0,1]$, given by $T_{b}(x)=b^{i}x-(b-1)$ for $x\in ((b-1)b^{-i},(b-1)b^{1-i}]$ and $i\in\mathbb{N}$. Using the argument in the last section with respect to $T_{b}$ instead of $T$ we obtain:
\begin{theorem}
For almost all $x=\langle n_{1},n_{2},n_{3},\dots\rangle_{b}\in(0,1]$ we have
\[ \lim_{i\to\infty} \sqrt[i]{n_{1}n_{2}\dots n_{i}}=\prod_{i=1}^\infty \sqrt[b^i]{i^{b-1}}=:\sigma_{b}.\]
\end{theorem}
We like to call $\sigma_{b}$ for $b>2$ a generalized Somos' constant. These constants are universal with respect to the base $b$ representation $\langle n_{1},n_{2},n_{3},\dots\rangle_{b}$. The generalization given here is slightly different from the generalization of Somos' constant studied in \cite{[SH]}, which is not related to universality.\footnote{They consider $\sqrt[b-1]{\sigma_{b}}$.} \\
We like to end the paper with a nice expression of generalized Somos' constants $\sigma_{b}$ using values of the generalized Euler-constant function
\[ \gamma(z)=\sum_{i=1}^{\infty}z^{i-1}\left(\frac{1}{i}-\log(\frac{i+1}{i})\right),\]
where $|z|\le 1$.
\begin{proposition} For all integers $b\ge 2$ we have
\[\sigma_{b}=\frac{b}{b-1}e^{-\gamma(1/b)/b}.\]
\end{proposition}
\proof We have
\[ \gamma(1/b)=\sum_{i=1}^{\infty}(b^{-i+1}/i-b^{-i+1}\log(i+1)+b^{-i+1}\log(i))\]
\[ =b( \sum_{i=1}^{\infty}b^{-i}/i-\sum_{i=1}^{\infty}b^{-i}\log(i+1)+\sum_{i=1}^{\infty}b^{-i}\log(i))\]
\[ =b(\log(b/(b-1))-b\sigma_{b}/(b-1)+\sigma_{b}/(b-1))=b\log(b/((b-1)\sigma_{b}))\]
using
\[ \sum_{i=1}^{\infty}b^{-i}/i=\log(b/(b-1)) \mbox{ and }\log(\sigma_{b})=\sum_{i=1}^{\infty}b^{-i}(b-1)\log(i).\]
Hence $e^{\gamma(1/b)}=(b/((b-1)\sigma_{b}))^{b}$ and $e^{\gamma(1/b)/b}=b/((b-1)\sigma_{b})$ given $\sigma_{b}=be^{-\gamma(1/b)/b}/(b-1)$.\qed~\\~\\
Estimates of $\gamma(1/b)$ can be found in \cite{[M]}.

\end{document}